\documentclass[10pt]{amsart}
\usepackage{amssymb, amscd, amsmath, amsthm, epsf, epsfig, latexsym,psfrag}

\usepackage{pb-diagram}
\def\zz{{\bf Z}}
 
\def\qq{{\bf Q}}
\def\cc{{\bf C}}
\def\rr{{\bf R}}
\def\co{\colon\thinspace}
\def\cs{\mathop{\#}}

\def\calg{\mathcal{G}}
\def\calc{\mathcal{C}}
\def\cala{\mathcal{A}}

\def\calj{\mathcal{J}}

\def\co{\colon}

\newcommand\Hom{\operatorname{Hom}}

\newcommand{\disc}{\operatorname{disc}}

\newtheorem{theorem}{Theorem}[section]
\newtheorem{lemma}[theorem]{Lemma}
\newtheorem{corollary}[theorem]{Corollary}
\newtheorem{prop}[theorem]{Proposition}

\theoremstyle{definition}
\newtheorem{definition}[theorem]{Definition}
\def\co{\colon\thinspace}

\numberwithin{equation}{section}

\begin{document}

\title[Non-slice    linear combinations of algebraic knots]{Non-slice    linear combinations of algebraic knots}
\author{Matthew Hedden}
\author{Paul Kirk}
\author{Charles Livingston} 
 
\thanks{   }
 \address{Matthew Hedden: Department of Mathematics, Michigan State University, East Lansing, MI 48824 \vskip.05in}
\address{Paul Kirk, Charles Livingston: Department of Mathematics, Indiana University, Bloomington, IN 47405 \vskip.1in}

\email{hedden@msu.edu}
\email{pkirk@indiana.edu}
\email{livingst@indiana.edu}
 


\begin{abstract} We show that the subgroup of the knot concordance group generated by links of isolated complex  singularities  intersects the subgroup of algebraically slice knots in an infinite rank subgroup.    \end{abstract}

\maketitle


\section{Introduction}

A long standing question~\cite{Rudolph1976} asks whether the set of algebraic knots, those that are links of isolated singularities of complex curves, is linearly independent   in the knot concordance group.   Following Rudolph's initial work~\cite{Rudolph1976}, Litherland~\cite{Lith} used signature functions to prove that the subset consisting of positive torus knots is independent.   However, it was then shown in~\cite{LM} by means of an example that invariants of the algebraic concordance group, such as the signature functions, are insufficient to prove the independence of algebraic knots. Somewhat later, Miyazaki~\cite{Miyazaki94} showed that the particular example in~\cite{LM} was not {\em ribbon}, and Rudolph~\cite{Rudolph2002} then used this result to conclude that ribbon knots do not behave well with respect to a certain operation called {\em plumbing}.     Here we develop tools  based on Casson-Gordon theory~\cite{CG} to prove the independence of a large family of algebraic knots, a family which includes the particular  example first found in~\cite{LM}.

To be more specific, let $\calc$ denote the group of (topologically locally flat) concordance classes of knots in $S^3$, and let $\calg$ denote the algebraic concordance group.  In~\cite{levine},  Levine constructed a surjection $\calc\to\calg$ whose counterpart in higher dimensions is an isomorphism.
The main result of this article is the following.

\medskip

\noindent{\bf Theorem 1.}  {\em Let $\cala$ denote the subgroup of the knot concordance group $\calc$ generated by algebraic knots: connected  links of isolated singularities in $\cc^2$. Then the intersection of $\cala$ with the kernel of $\calc\to \calg$ contains an infinitely generated free abelian subgroup.}
 \medskip

We use the following notation: given a knot $K$ and a pair of relatively prime integers $p$ and $q$, $K_{p,q}$ denotes the oriented $(p,q)$--cable of $K$.   Thus $K_{p,q}$ represents $p$ times the generator of the first homology of  the tubular neighborhood of $K$.   In the special case that $K$ is the unknot $U$, so that $K_{p,q}$ is a torus knot, we use the standard notation of $T_{p,q}$ rather than $U_{p,q}$.  The notation can be iterated; for instance, $K_{p,q;r,s}$ denotes the $(r,s)$--cable of the $(p,q)$--cable of $K$.

 An {\em algebraic knot} is, by definition, the connected link of an isolated singularity of a polynomial map $f \co \cc^2\to \cc$. A knot is isotopic to  an algebraic knot if and only if it is an iterated torus knot $T_{p_1,q_1;\cdots; p_n,q_n}$ with indices satisfying $p_i, q_i>0$ and $q_{i+1}>p_iq_i p_{i+1}$ (see, for instance,~\cite{EN1985}).

 Our results concern  $(2,k)$--cables of knots.  In particular, we resolve an old question of whether a particular linear combination of $(2,k)$--cables is slice; this combination is the simplest algebraically slice knot in the span of the algebraic knots,~\cite{LM}:
 
 \medskip

 \noindent{\bf Theorem 2.}  {\em The linear combination of algebraic knots 
 $$ T_{2,3;2, 13} \cs T_{2,15} \cs -T_{2,3;2,15} \cs  -T_{2,13} $$
 is algebraically slice but has infinite order in $\calc$.}
 
 \medskip
 
Theorems 1 and 2 are consequences of the following result, which establishes the linear independence of an infinite collection of algebraic knots.

\medskip

\noindent{\bf Theorem 3.} {\em For appropriately chosen integers $q_i$, the set of algebraic knots 
 $$\{ T_{2,q_i}, T_{2,3;2,q_i}\}_{i=1}^\infty$$
 form a basis of a free abelian subgroup of  the concordance group $\calc$. This subgroup intersects the kernel of $\calc\to\calg$ in a free abelian subgroup, with basis given by  the following set of algebraically slice knots
 $$\{ T_{2,3;2, q_{n}} \cs T_{2,q_{1}}\cs -T_{2,3;2,q_{1}}\cs  -T_{2,q_{n}}\}_{n=2}^{\infty} .$$
 }
 
 \medskip
 
  Our arguments apply more generally, for example  to cables of knots other than  the trefoil $T_{2,3}$.  
  We refer the reader to the body of the paper for details.

\medskip
 
 The methods we use are those introduced by Casson-Gordon in~\cite{CG}. A novel feature of our approach is the essential interplay between signature and discriminant invariants on the Witt group of Hermitian forms over $\cc(t)$.  Casson-Gordon signature invariants, which are $\zz$--valued and hence more effective in identifying elements of infinite order, are often intractable to compute.  Discriminant invariants are computable algorithmically~\cite{KL1}, but because they take values in a group that is $\zz/2\zz$--torsion, they are less effective in determining linear independence.   By combining the two types of invariants  we are able to apply the power of signatures while bypassing the need to explicitly compute their values; we can also avoid most of the typically messy work of analyzing metabolizers which plagues many discriminant arguments used to show that certain knots have infinite order in $\calc$.
 
 \medskip
 
 We finish this introduction by giving some background to place our results in context.

Given an oriented knot $K$, let $-K$ denote the mirror image of $K$ with its orientation reversed. Two oriented knots, $K_1$ and $ K_2$, are called {\em concordant} if the connected sum $K_1\cs  -K_2 $ bounds a locally flat embedded disk in $S^3$.  The set of concordance classes of knots forms an abelian group $\calc$  with operation induced by connected sum.  Knots which represent the zero element in $\calc$ are called {\em slice}.

From this point forward  we do not distinguish in our notation between a knot and its concordance class.  In particular, we write $K_1 - K_2$ for the connected sum $K_1 \cs - K_2$.  We will also write $-K_{p,q}$ for $ -(K _{p,q})$, both of which equal $(-K)_{p,-q}$.

 Fox and Milnor observed~\cite{FoxMilnor} that if two knots are concordant, then the product of their Alexander polynomials is a norm    in $\zz[t^\pm]$: that is, 
 $\Delta_{K_1}(t)\Delta_{K_2}(t)=  f(t)f(t^{-1})$ for some polynomial $f(t)\in  \zz[t^\pm]$. (Recall that the Alexander polynomial is defined up to multiplication by $\pm t^{k}$.) 
 
  An early result of Seifert~\cite{Seifert1950}  
  (see~\cite[Theorem 6.15]{Lickorish} for a recent reference) 
  states that the Alexander polynomial of a satellite knot is determined by the Alexander polynomials of the knots involved in the construction, together with an integer called the winding number.  In the case of cables, the formula is given by   
   $$\Delta_{K_{p,q}}(t) = \Delta_{T_{p,q}}(t) \Delta_K(t^p), \ \ \text{where}\ \   \Delta_{T_{p,q}}(t) =\frac{(t^{pq} -1)(t-1)}{(t^p-1)(t^q-1)}.$$
   For a connected sum, the Alexander polynomial is simply the product of the Alexander polynomials of the constituent knots. 
 A bit of calculation using these facts   shows that distinct  algebraic knots are not concordant. The  Levine--Tristram signatures  of a knot~\cite{levine,  Lickorish, mi, tr}, which define  integer-valued homomorphisms on $\calc$, can be used to further show that algebraic knots have infinite order in $\calc$.  

These observations might  lead one to conjecture that the set of algebraic knots forms a basis for an infinitely generated free abelian subgroup $\cala\subset \calc$.  A first line of attack to this question, as taken in~\cite{Litherland1979} and~\cite{LM}, is to consider the {\em algebraic concordance group} $\calg$ and to determine the image of the composite $\cala\subset\calc\to\calg$. For the purposes of this article the precise definition of $\calg$ is not needed  and it will suffice to say that $\calg$ is the group generated by Seifert forms of knots, modulo Seifert forms of slice knots.  The relevant facts surrounding $\calg$ are: 

 \begin{enumerate}
\item There is a surjection $\calc\to\calg$~\cite{levine}.
\item The algebraic concordance class of a knot $K$ is  determined by its Blanchfield (torsion) form~\cite{Trotter1973}:
$$Bl_K \co  H_1(S^3-K;\zz[t^\pm])\times H_1(S^3-K;\zz[t^\pm])\to\frac{\qq(t)}{\zz[t^\pm]}.$$
 
\end{enumerate}

With respect to the interplay of cabling and algebraic concordance, the   formula~\cite{Ke, LivingstonMelvin1985} 
 $$Bl_{K_{p,q}}(t) = Bl_K(t^p) \oplus Bl_{T_{p,q}}(t),$$
based on a Mayer-Vietoris argument,     gives a quick method to determine  if certain linear combinations of cable knots are  {\em algebraically slice}, that is lie in the kernel of $\calc\to \calg$.

As observed in~\cite{LM} (see Lemma \ref{summarylemma}  below), this formula implies that the knot in Theorem 2,   $ T_{2,3;2, 13} + T_{2,15}-T_{2,3;2,15} -T_{2,13}$, is algebraically slice. It represents the simplest example of a knot in the kernel of the composite $\cala\subset\calc\to\calg$.   Showing this knot is not slice has remained open until now, although Miyazaki proved it is not ribbon~\cite{Miyazaki94}.  

The reader will have noticed that the term ``algebraic'' has two different meanings in this paper; on the one hand it describes a class of knots defined as links of isolated singularities, and on the other it describes a certain quotient of the knot concordance group.  Algebraic knots are iterated cables,  and we will typically work with general cables, so this should cause no confusion.

\subsection{Comparison with smooth techniques}
Progress in identifying the structure of algebraic knots in the setting of smooth concordance has been achieved largely through analytic means or the deep combinatorial approach stemming from Khovanov homology theory.  This is most notable in the   solutions to the Milnor conjecture  and the proof that the smooth 4--ball genus of a torus knot is realized by an algebraic curve~\cite{km, os, ra}.  Highlighting the necessity of smooth techniques in studying algebraic knots, Rudolph~\cite{Rudolph93} observed that the Milnor conjecture is false in the topological locally flat category. 

Thus, it comes as a surprise that Casson-Gordon methods apply so effectively in the present setting, having the further advantage that we can establish the independence of these knots in the topological concordance group. 
Nonetheless,  it would be interesting to know the extent  to which the array of existing smooth concordance invariants can be used to address the question of independence of algebraic knots. 
We should point out, however, that the Ozsv{\'a}th-Szab{\'o}~\cite{os} and Rasmussen~\cite{ra} concordance invariants, $\tau$ and $s$,  coming from  knot Floer homology and Khovanov homology, respectively, contain no information for the knots at hand.  We make this precise in Proposition \ref{prop:stau}, which shows that both invariants vanish for the family above and its obvious generalization to positively iterated torus knots.  

Despite the failure of $s$ and $\tau$, it seems likely that grading information from the Floer homology of branched covers (in the form of the Fr{\o}yshov invariant~\cite{Froyshov96}, Ozsv{\'a}th-Szab{\'o} correction terms~\cite{AbsGrad}, or Chern-Simons invariant of $SU(2)$ representations~\cite{MR1051101, MR1047138}) could be useful in our pursuit.   However,  extensive computations of such invariants in the first two cases is difficult, and  computing Chern-Simons invariants of covers in the spirit of Fintushel-Stern~\cite[Theorem 5.1]{MR1051101} and Furuta~\cite[Theorem 2.1]{MR1047138}    have failed to determine if any member of the family of knots in the present article are slice.

\bigskip

\noindent{\bf Acknowledgment:} Conversations with Tom Mrowka about the ribbon-slice conjecture led naturally to the investigations of this paper.  The relationship between this ribbon-slice problem and the concordance independence of algebraic knots is  discussed in the second appendix.


\section{Two--fold branched covers and characters}
 As mentioned in the introduction,  we write $K_{p,q}$ for the $(p,q)$--cable of $K$  and $-K_{p,q}$ for $-(K_{p,q})$, which, by a simple orientation argument, equals $(-K)_{p,-q}$. Cabling the concordance shows that the concordance class of $K_{p,q}$ depends only on the concordance class of $K$, in the sense that if $K$ and $K'$ are concordant, then $K_{p,q}$ and $K'_{p,q}$ are concordant.  These observations, along with our earlier description of the Blanchfield pairing of a cable knot, yields the following general statement, implicit in~\cite{LM}.

\begin{lemma}\label{summarylemma} For any knot K, 
 $$K_{p,q_1}   + T_{p, q_2}  - K_{p, q_2}   - T_{p, q_1} $$ 
is an algebraically slice knot    and is a slice knot when $K$ is slice.\qed
\end{lemma}

\begin{figure}
\psfrag{K}{$q$}
\psfrag{U}{$U$}
\begin{center}
\includegraphics[height=110pt]{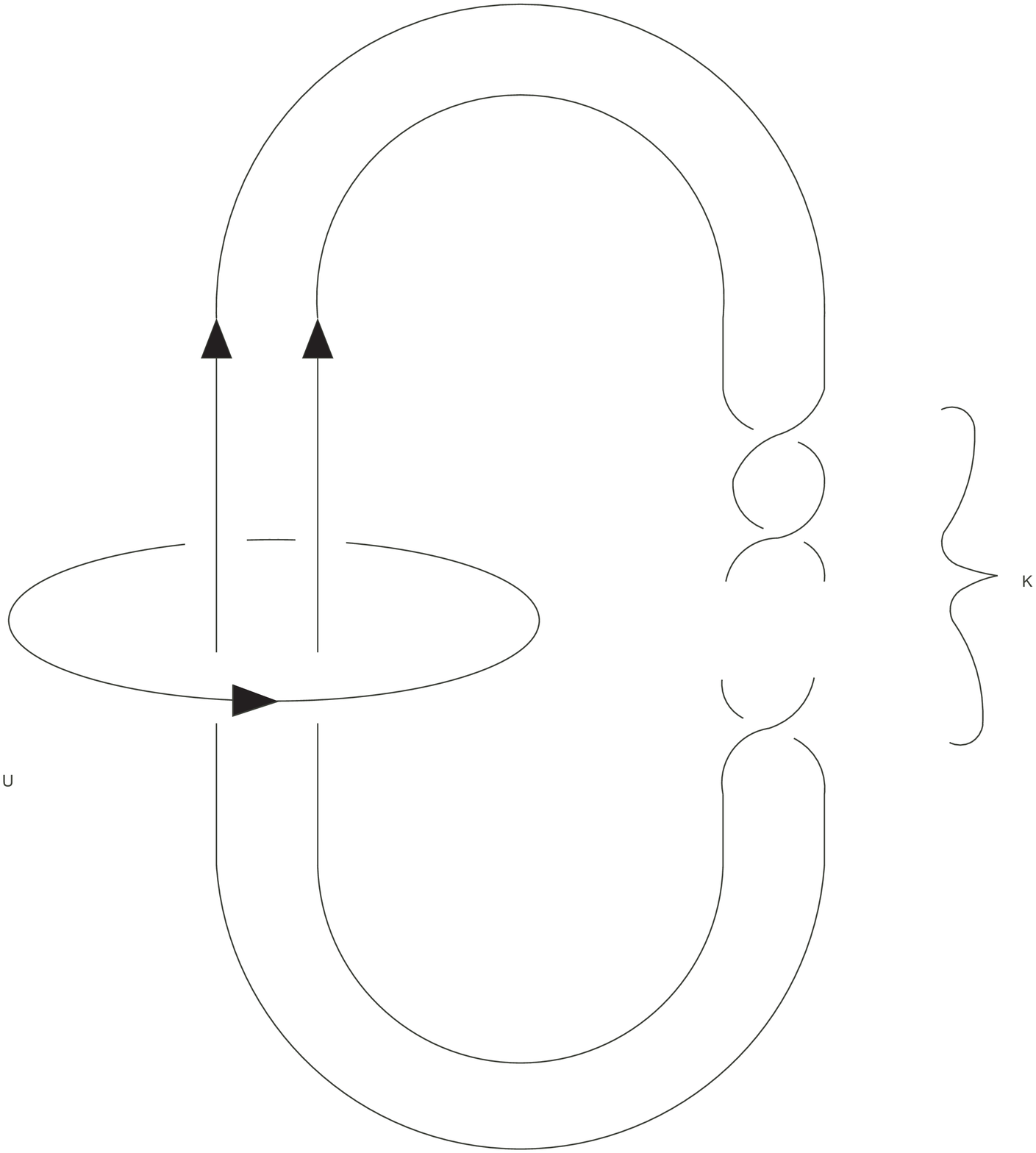}
\caption{\label{fig:Lq}$L_q$}
\end{center}
\end{figure}

 We now turn our focus to 2--stranded cables; knots of the form $K_{2,q}$. 
A useful depiction of $K_{2,q}$ is the following. Figure \ref{fig:Lq} shows a 2--component link $L_q$ with one component a $(2,q)$--torus knot and the other component an unknot labeled $U$.   If $K$ is a knot in $S^3$  and $q$ is an odd integer, then $K_{2,q}$ is obtained by removing a neighborhood of $U$ and replacing it by the complement of a tubular neighborhood of $K$ in such a way that the meridian-longitude pairs of $U$ and $K$ are interchanged.  

For a knot $K$, denote by   $M^3(K)$ the   2--fold branched cover of $S^3$  branched over $K$. Let $\widetilde{K}$ denote the lift of $K$ to $M^3(K)$ and let $M^3_0(K)$ denote the result of $0$--surgery on $M^3(K)$ along $\widetilde{K}$; that is, the surgery whose framing comes from a lift of the longitude of $K$.  Note that $M^3_0(K)$ is the 2--fold cyclic cover of 0--surgery on $K$ in $S^3$.

\begin{figure}
\psfrag{U1}{$\tilde{U}_1$}
\psfrag{U2}{$\tilde{U}_2$}
\psfrag{K}{$q$}
\begin{center}
\includegraphics[height=110pt]{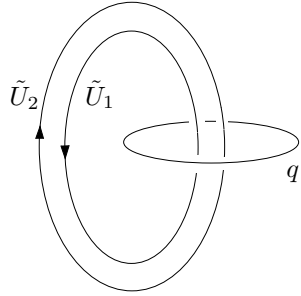}
\caption{\label{fig:BkU}$M^3(T_{2,q})$.  The branched double cover of the $(2,q)$--torus knot is the lens space $L^3(q,1)$, obtained by performing $q$--surgery on an unknot.  The unknotted component of $L_q$ lifts to $\widetilde{U}_1\cup \widetilde{U}_2$.}
\end{center}
\end{figure}

The 2--fold branched cover $M^3(T_{2,q})$  is  the lens space $L^3(q,1)$.  
Since $U$ links the $(2,q)$--torus knot twice in $L_q$, the preimage of $U$ in this 2--fold branched cover consists of two curves, $\widetilde{U}_1$ and   $\widetilde{U}_2$.  One way to understand this is to take a 3--ball in $S^3$ which meets  the $(2,q)$--torus knot in two unknotted arcs and contains $U$ in its interior.  Then the preimage of  this 3--ball in $M^3(T_{2,q})$ is a solid torus, as is the preimage of its complement in $S^3$.   The curves $\widetilde{U}_1$ and   $\widetilde{U}_2$ are each circles parallel to the core of this solid torus but oppositely oriented.  In particular, $\widetilde{U}_1$ and   $\widetilde{U}_2$ are isotopic  as unoriented curves   in $M^3(T_{2,q})$. Figure \ref{fig:BkU} depicts the situation.

To obtain $M^3(K_{2,q})$, we replace the solid torus neighborhood of $\widetilde{U}_1$ and   $\widetilde{U}_2$     with copies  of the complement of $K$ in $S^3$, interchanging the meridian-longitude pairs.  The preimage $\widetilde{K} \subset M^3(T_{2,q})$ is not drawn in Figure \ref{fig:BkU}.

We will need notation for certain curves in $M^3(K_{2,q})$.  
Denote by $\mu_U$ and $\lambda_U$ the meridian and longitude of the unknotted component $U\subset L_q$, with its orientation as in Figure \ref{fig:Lq}.      From the perspective of $K_{2,q}$ as a satellite knot, these are the longitude and meridian, respectively, for the  companion, $K$. 
Denote by $\tilde{\mu}_{i}$ and   $\tilde{\lambda}_{i}$ the meridian  and longitude  in the boundary of the tubular neighborhood of $\widetilde{U}_i\subset M^3(T_{2,q})$, $i=1,2$ in the surgery diagram given in Figure \ref{fig:BkU}.    

The notation is somewhat ambiguous since there is no preferred choice of lift of $U$, but in any case we will choose the same ordering when comparing $M^3(K_{2,q})$ and $M^3(T_{2,q})$.  Note that  $\tilde{\mu}_1$ and $\tilde{\mu}_2$ vanish in $H_1(M^3(K_{2,q}))=\zz/q\zz$,  $\tilde{\lambda}_1$ generates $H_1(M^3(K_{2,q}))$, and $\tilde{\lambda}_2=-\tilde{\lambda}_1$ in $H_1(M^3(K_{2,q}))$.

 Denote by $\tilde{\mu}$ the preimage of the meridian to $K_{2,q}$ and by   $\tilde{\lambda}$ a component of the preimage of the longitude  of $K_{2,q}$ in $M^3(K_{2,q})$.   In particular, $\tilde{\mu}$ and $\tilde{\lambda}$ are nullhomologous in $M^3(K_{2,q})$, since a Seifert surface for $K_{2,q} \subset S^3$ lifts to $M^3(K_{2,q})$. 

The linking form of $M^3(K_{2,q})$ is given by the $1\times 1$ matrix $(\tfrac{1}{q})$; in fact  
$$lk(\tilde{\lambda}_1,\tilde{\lambda}_1)=1/q=lk(\tilde{\lambda}_2,\tilde{\lambda}_2).$$

Let  $p$ be an odd prime and let $C_p$ denote the cyclic group of order $p$. This group can be identified with the group of   $p$--roots of unity: $C_p = \{\zeta_p^a\}\subset \cc$, where $\zeta_p=e^{2\pi i/p}$.  If $p$ and $q$ are relatively prime then every character $\chi \co H_1(M^3(K_{2,q}))\to C_p$ is trivial. 
On the other hand, if $p$ divides $q$, then 
the set of all characters $\chi \co H_1(M^3(K_{2,q})) \to C_p$ form a cyclic group isomorphic to $C_p$. We can fix an isomorphism $\Hom(H_1(M^3(K_{2,q})),C_p)\cong C_p$ as follows:  let  
 $$\chi_1 \co H_1(M^3(K_{2,q}))\to C_p$$ denote the character which takes $\tilde{\lambda}_1$ to  $\zeta_p$. Then any other character is obtained by post-composing $\chi_1$ with the homomorphism $C_p\to C_p$ of the form $\zeta_p^i \mapsto \zeta_p^{ia}$ for some integer $a$.   We denote this composite as $\chi_a \co H_1(M^3(K_{2,q}))\to C_p$.
Notice that $a$ is well defined modulo $p$, and,  although the definition of $a$ depends on a choice of ordering of the two lifts $\widetilde{U}_1, \widetilde{U}_2$,  the unordered pair $\{a, -a\}$ is independent of this choice.

Then 
$$\chi_a(\tilde{\lambda}_1)=\zeta_p^a, \ \chi_a(\tilde{\lambda}_2)=\zeta_p^{-a}, \ \chi_a(\tilde{\mu}_1)=1, \ \chi_a(\tilde{\mu}_2)=1, \ \chi_a(\tilde{\mu})=1, 
\text{ and } \chi_a(\tilde{\lambda})=1.$$

  The  unbranched   2--fold cover $M^3(K_{2,q})-\widetilde{K}_{2,q}\to  S^3-K_{2,q}$ induces a homomorphism $H_1(M^3(K_{2,q}) -\widetilde{K}_{2,q})\to H_1(S^3-K_{2,q})=\zz$ with image $2\zz$. Dividing by two defines a surjection $\epsilon \co H_1(M^3(K_{2,q}) -\widetilde{K}_{2,q})\to \zz$. Writing $\zz=\langle t\rangle$ multiplicatively, we have
  $$\epsilon(\mu_i)=1, \epsilon(\tilde{\lambda}_i)=t, \epsilon(\tilde{\mu})=t,\text{ and } \epsilon(\tilde{\lambda})=1.$$  To see that $\epsilon(\tilde{\lambda}_1)=t=\epsilon(\tilde{\lambda}_2)$, notice that $\tilde{\lambda}_i$ is sent to $\lambda_U$ in $S^3-K_{2,q}$, which links $K_{2,q}$ twice; dividing by two  yields one.

  Recall that $M^3_0(K_{2,q})$ denotes the closed 3--manifold obtained by performing $0$--surgery on  $\widetilde{K}_{2,q} \subset M^3(K_{2,q})$; that is, the surgery corresponding to the framing induced by a lift of a longitude of $K_{2,q}$ to $M^3(K_{2,q})$. 
    Since $\chi_a(\tilde{\lambda})=\zeta_p^0=1$ and $\epsilon(\tilde{\lambda})=t^0=1$,  both $\chi_a$ and $\epsilon$ uniquely  extend to homomorphisms on $H_1(M^3_0(K_{2,q}))$.  We can view their product as a homomorphism to the multiplicative group of non-zero elements of the field of rational functions, $\cc(t)$:
  $$\chi_a\times \epsilon \co H_1(M^3_0(K_{2,q}))\to \cc(t)^\times.$$ Each homology class is sent to an element of the form $\zeta_p^b t^c$. 
  
 We summarize these facts  in the following lemma.   
 \begin{lemma}   Let $K$ be a knot in $S^3$, $K_{2,q}$ its  $(2,q)$--cable, and $T_{2,q}$ the $(2,q)$--torus knot.     Let $M^3(K_{2,q})$ denote the 2--fold branched cover of $S^3$ branched over $K_{2,q}$, and let $M^3_0(K_{2,q})$ denote the manifold obtained from 0--surgery on $M^3(K_{2,q})$ along the preimage of the branch set.  Choose an odd prime $p$ and let $\zeta_p=e^{2\pi i/p}$.  Then
 \begin{enumerate}
 \item $M^3(K_{2,q})$ is obtained from $M^3(T_{2,q})=L^3(q,1)$ by removing   neighborhoods of the two preimages  $\widetilde{U}_1$, $\widetilde{U}_2$ of $U$ and gluing in two copies of $S^3-nbhd(K)$, so that the meridian-longitude pairs of $\widetilde{U}_i$ and $K$ are interchanged.
 
 \item   $H_1(M^3(K_{2,q}))=\zz/q\zz$, generated by $\tilde{\lambda}_1$, and $lk(\tilde{\lambda}_1,\tilde{\lambda}_1)=1/q$.
\item To any character $\chi \co H_1(M^3(K_{2,q}))\to C_p$ one can associate the integer $a$ modulo $p$ by the condition $\chi(\tilde{\lambda}_1)=\zeta_p^a$.  This character is denoted $\chi_a$. In particular, this sets up a 1-1 correspondence between $C_p$--valued characters on $H_1(M^3(K_{2,q}))$ and  on $H_1(M^3(T_{2,q}))$.
\item The character $\chi_a$ uniquely determines a character (also denoted $\chi_a$) on $H_1(M^3_0(K_{2,q}))$.
\item There is a surjection $\epsilon \co H_1(M^3_0(K_{2,q}))\to \zz =  \langle t\rangle$ satisfying
$\epsilon(\mu_i)=1, \epsilon(\tilde{\lambda}_i)=t, \epsilon(\tilde{\mu})=t,\text{ and } \epsilon(\tilde{\lambda})=1.$

\end{enumerate}
\qed

\end{lemma}


\section{Casson-Gordon invariants}

Let $\calj \co \cc(t)\to \cc(t)$ denote the involution $\calj(f(t))=\bar{f}(t^{-1})$; specifically,
$$ \calj \big(\frac{\sum a_i t^i}{\sum b_j t^j} \big)= \frac{\sum  \overline{a}_i t^{-i}}{\sum \overline{b}_j t^{-j} },$$ where $\overline{a}_i$ denotes complex conjugation. 
 We let $W(\cc(t), \calj)$ denote the corresponding Witt group of non-singular $\calj$--Hermitian forms.  This Witt group is discussed in more detail in Appendix~\ref{appwitt}.  In brief, two forms $I_1$ and $I_2$ are equivalent if the sum $I_1\oplus  -I_2 $ is metabolic; that is, if it contains a half-dimensional   subspace on which the form vanishes.    The set of equivalence  classes of forms constitute the Witt group,  with operation induced by direct sums.

 For each choice of $K$ and $\chi \co H_1(M^3(K_{2,q}))\to C_p$ with $p$ an odd prime, the {\em Casson-Gordon invariant  of $(K_{2,q}, \chi)$}, 
 $$\tau (K_{2,q}, \chi)\in W(\cc(t), \calj)\otimes \zz_{(2)},$$
   is defined as follows~\cite{CG}.   (Here $\zz_{(2)}$ is  $\zz$ {\it localized at} 2: the set of rational numbers with odd denominator.)   Elementary bordism theory shows that  $p \cdot (M^3_0(K_{2,q}), \chi\times \epsilon)$ bounds: say $p \cdot(M^3_0(K_{2,q}), \chi\times \epsilon)= \partial(Y^4, \rho)$.  Then $Y^4$ has a non-singular $\cc(t)$--valued intersection form $I(Y^4,\rho)\in W(\cc(t),\calj)$ defined using the cup product on middle degree cohomology, with  local coefficients determined by the homomorphism $\rho \co \pi_1(Y^4)\to \cc(t)^\times$.  On the other hand, $Y^4$ also has  its ordinary intersection form $I(Y^4)\in \text{ Image } \{W(\qq)\to W(\cc(t), \calj)\}$.  The  Casson-Gordon invariant is defined to be
$$\tau (K_{2,q}, \chi)=\tfrac{1}{p}(I(Y^4,\rho)-I(Y^4)).$$ Since $p$ is odd,  $\frac{1}{p} \in \zz_{(2)}$.

The  correspondence between characters on $M^3(T_{2,q})$ and $M^3(K_{2,q})$ described in the previous section permits us to unambiguously define the difference $\tau(K_{2,q} ,\chi)-\tau(T_{2,q},\chi)$.   A formula for this difference was established  (in much greater generality) by Litherland in his influential article~\cite{Lith}.  (See also Gilmer~\cite{Gilmer} for related work and applications of this approach.) 
Using this result, one can compute the difference of  Casson-Gordon invariants for different choices of $K$.   The answer is given in terms of an abelian invariant, $\alpha_K$, which we define next.

Let $S_0^3(K)$ denote the 3--manifold obtained by $0$--surgery on the knot $K\subset S^3$. The orientation of $S^3$ and $K$ determine an isomorphism $\delta \co H_1(S_0^3(K))\to \zz=\langle x\rangle$.  There is a   4--manifold $X^4$ and $\bar{\delta} \co \pi_1(X^4)\to \langle x\rangle$ so that $\partial(X^4,\bar{\delta})=(S_0^3(K),\delta)$.  Then $X^4$ has a $ \qq[x^\pm]$-equivariant intersection form $I(X^4,\bar{\delta})$ and an ordinary integer-valued intersection form $I(X)$.  The concordance invariant, $\alpha_K$, is defined to be the difference of these forms in the Witt group of $\qq(x)$:
$$\alpha_K =  I(X^4,\bar{\delta}) - I(X^4) \in W(\qq(x), \calj).$$

The class $\alpha_K\in W(\qq(x), \calj)$ is determined by the algebraic concordance class of $K$; that is, the image of $K$ under the map $\calc\to\calg$.   Given a unit complex number $\omega$, the {\em Levine--Tristram $\omega$--signature of $K$} is defined to be the signature of the complex Hermitian matrix obtained by substituting $x=\omega$ into a matrix representative of $\alpha_K$.

More generally, if $\omega $ is a unit complex number, the map $x \to \omega t$ induces a map $W(\qq(x), \calj) \to W(\cc(t), \calj)$.    We define 
$\alpha_K(\omega t)$ to be the image of $\alpha_K$ under this map.

Litherland's theorem~\cite[Corollary 2]{Lith}, proven by a delicate Mayer--Vietoris argument,   implies the following.

\begin{prop} \label{prop2.1} Given $K, q, \chi$, and $p$ as above, then
$$\tau(K_{2,q},\chi_a)-\tau(T_{2,q},\chi_a)= \alpha_K(\zeta_p^a t) + \alpha_K(\zeta_p^{-a} t) 
$$
in $W(\cc(t))\otimes \zz_{(2)}$.\qed
\end{prop}

  Notice   that 
$ \alpha_{K}(\zeta_p^{a} t) +
 \alpha_{K}(\zeta_p^{-a} t)$ is unchanged by replacing $a$ by $-a$.   
Moreover, for any knot $K$ and character $\chi$, (writing $\chi$ additively in this formula for simplicity)
 \begin{equation}\label{eq4.4}
  \tau(K,\chi)=\tau(K, -\chi).\end{equation}
This is because the  2--fold covering  transformation is an orientation-preserving diffeomorphism which preserves the orientation of the branch set  but induces $-1$ on the first homology of the branched cover.  Precomposing $\chi$ with this diffeomorphism yields $-\chi$. Hence  $\tau(K_{2,q}, \chi_{a  })= \tau(K_{2,q}, \chi_{-a })$.  

  In particular, to a character $\chi \co H_1(M^3(K_{2,q}))\to C_p$ we can unambiguously assign  $a\in\{0,1,2,\cdots, \tfrac{p-1}{2}\}$ by evaluating 
 $\chi$ on one of $\tilde{\lambda}_1$ or $\tilde{\lambda}_2$ and replacing $\zeta_p^a$ by $\zeta_p^{p-a}=\zeta_p^{-a}$ if necessary.
 The number $a$   determines $\chi$   up to sign, but it completely determines $\tau(K_{2,q},\chi_a)$ and 
 $ \alpha_{K}(\zeta_p^{a} t) +
 \alpha_{K}(\zeta_p^{-a} t)$.  This also resolves the ambiguity introduced earlier in choosing an order of the lifts of $U$, since if $\chi(\tilde{\lambda}_1)=a$, then $\chi(\tilde{\lambda}_2)=-a$.  Notice   that $\chi$ is trivial if and only if $a=0$.
 
 \vskip.1in
 
We conclude this section with a lemma describing the role of orientation on the value of $\tau$ and $\alpha$.
\begin{lemma} $\alpha_K = -\alpha_{-K}$ and $\tau(K, \chi) = -\tau(-K, \chi)$.
\end{lemma}
\begin{proof} If we consider a representative of $K$ to be an embedded $S^1$ in $S^3$, then $-K$ is represented by the same $S^1 $ in $S^3$, but with the orientation of $S^3$ (and of $S^1$) reversed.  Hence, there is a natural orientation-reversing homeomorphism from $M^3(K)$ to $M^3(-K)$.  This permits us to formally make sense of the statement of the lemma; characters on the covers of $K$ and   $-K$ can be identified via this homeomorphism.  

Given this, the only difference between the computation of the Witt class invariants of $K$ and $-K$ are that the relevant 4--manifolds have their orientations reversed.  This has the effect of changing the signs of the intersection forms.  
\end{proof}


\section{Linear combinations and slicing}\label{sectionli}

  Casson-Gordon invariants are used to obstruct sliceness of knots. The main result of~\cite{CG} implies that 
if a knot $K$ is slice, then there exists a metabolizer $V\subset H_1(M^3(K))$ for the linking form on $M^3(K)$ (as earlier, $M^3(K)$ denotes the 2--fold branched cover of $S^3$ over $K$), so that $\tau(K,\chi)=0$ for every character $\chi \co H_1(M^3(K )) \to \cc^\times$ that factors through $\zz/p \zz  $ and   vanishes on $V$.   (Recall, a {\em metabolizer} is a subgroup $V\subset H_1(M^3(K))$ on which the linking form vanishes  and for which the order of $V$ is the square root of the order of $H_1(M^3(K))$.)
  
 If $p$ is a prime dividing the order of 
 $H_1(M^3(K))$, then given any metabolizer $V$,  one can find a {\em non-trivial} $C_p$--valued character which vanishes on $V$, since $H_1(M^3(K))/V$ necessarily has  non-trivial $p$--torsion.  Therefore, given any knot $K$ and  a prime $p$  dividing the order of $H_1(M^3(K))$,  if $\tau(K,\chi)\ne 0$ for all non-trivial $C_p$ characters $\chi$,   $K$ is not slice.

Suppose now that  sequences $K^i$ of knots and  $q_i$ of relatively prime odd integers,  $i=1,2,\cdots$, are given.    Although our techniques apply more generally, for our applications we can assume that the $q_{2i-1}$ are prime, so henceforth do.

 As explained above, 
the linear combination of   cables 
\begin{equation}\label{eq3.1}
J_i= K^i_{ 2, q_{2i-1}} + T_{ 2, q_{2i}}   -K^i_{2,  q_{2i}}  -T_{ 2, q_{2i-1}}
\end{equation} 
is algebraically slice. 

\begin{lemma} \label{lem3.1}  If $\theta$ denotes the trivial character on $M^3(J_i)$, the 2--fold branched cover of $S^3 $ branched over $J_i$,   then $
\tau(J_i,\theta)=0$.
\end{lemma} 

\begin{proof} Let $\theta$ denote the trivial character on the first homology of the 2--fold branched cover of any knot. Applying Proposition \ref{prop2.1} and using the fact that the Casson-Gordon invariants are additive with respect to connected sums of knots (see, for instance,~\cite[Corollary 1]{Lith} or~\cite{Gilmer})  one computes
\begin{eqnarray*}\tau(J_i,\theta)&=&\tau(T_{2,  k_{2i-1}} ,\theta) +\tau(T_{ 2, k_{2i}},\theta)+ \tau(-T_{2,  k_{2i}},\theta)+ \tau(-T_{ 2,k_{2i-1}},\theta)\\
&&\quad + 2\alpha_{K^i}(t)+2\alpha_{-K^i}(t).
\end{eqnarray*}
But $\alpha_{K^i}(t)=-\alpha_{-K^i}(t)$ and, since $\theta$ is trivial, $\tau(T_{2,k},\theta)=-\tau(-T_{2,k},\theta)$. The lemma follows.
\end{proof}

Consider an algebraically slice linear combination  
\begin{equation}\label{eq3.2}J=\sum_{i=1}^N n_iJ_i.\end{equation}
The 2--fold branched cover $M^3(J)$ of $S^3$ branched over $J$ is the connected sum of the (oriented) branched covers of the 
constituent knots in $J$.  Hence
$$M^3(J)= \cs_{i=1}^N  n_i  \big( M^3(K^i_{2,  q_{2i-1}}) \cs M^3(T_{ 2, q_{2i}})\cs M^3(-K^i_{2,  q_{2i}})\cs M^3(-T_{ 2, q_{2i-1}})\big).$$
Assume  that  $  n_1>0.$ Let
$\chi \co H_1(M^3(J))\to C_{q_1}$ be a character. Let  $\zeta_{q_1}=e^{2\pi i/{q_1}}$.
The assumption that $q_i$ is relatively prime to $q_1$ for $i>1$ implies that  $\chi$   vanishes on each summand in the connected sum, except possibly for some of the $M^3(K^1_{2,q_1}) $ and $M^3(-T_{2,q_1} )$ summands.   On these summands, $\chi$ determines integers $a_1,a_2,\cdots, a_{n_1}$ and $b_1,\cdots, b_{n_1}$ in $\{ 0,1,2,\cdots, \tfrac{q_1-1}{2}\}$ by restricting $\chi$ to the $M^3(K^1_{2,q_1}) $ and $M^3(-T_{2,q_1} )$ summands, respectively, and evaluating on the corresponding lifts $\tilde{\lambda}_1$ or $\tilde{\lambda}_2$  in each summand, as in the previous section.  Using Lemma~\ref{lem3.1}, one concludes
\begin{equation}\label{eq3.3}
 \tau(J,\chi)= 
n_1\big(\tau(-K^1_{2, q_2} ,\theta)+\tau(T_{2,q_2},\theta)  \big) +
\sum_{i=1}^{n_1}\big( \tau(K^1_{2,q_1},\chi_{a_i})+\tau( -T_{2,q_1},\chi_{b_i})\big)\end{equation}
 where $\chi_{a_i}$ denotes the restriction of $\chi$ to $H_1(M^3(K^1_{2,q_1})) $ and $\chi_{b_i}$ denotes the restriction of $\chi$ to $H_1(M^3(-T_{2,q_1}) )$.

 Proposition \ref{prop2.1} gives the two equations: 
 $$\tau(-K^1_{2,q_2},\theta)+\tau(T_{2, q_2},\theta)=  
 2\alpha_{-K^1}(t)
 $$
 $$\tau(K^1_{2,q_1},\chi_{a_i})=\tau(T_{2, q_1},\chi_{a_i})+ \alpha_{K^1}(\zeta_{q_1}^{a_i} t)+ \alpha_{K^1}(\zeta_{q_1}^{-a_i} t).$$

Substituting these equations in Equation (\ref{eq3.3}) shows that

\begin{equation}\label{eq4.3}
\tau(J,\chi)=2n_1 \alpha_{-K^1}(t)  +\sum_{i=1}^{n_1}\big(  \alpha_{K^1}(\zeta_{q_1}^{a_i} t) +
 \alpha_{K^1}(\zeta_{q_1}^{-a_i} t) + \tau(T_{2,k_1},\chi_{a_i})-\tau(T_{2,k_1},\chi_{b_i})
 \big). 
 \end{equation} 

Summarizing, we have:

 \begin{prop} \label{prop3.2} Let   $J_i$ and $J$ be the knots described in Equations (\ref{eq3.1}) and (\ref{eq3.2}). Assume that $n_1 > 0$ and that $p$ is an odd prime not dividing $q_i$ for $i > 1$, and $\chi \co H_1(M^3(J))\to C_{p}$  a character, determining integers $a_i, b_i\in\{0,1,\cdots, \tfrac{p-1}{2}\}$ as described above.    Then

  \begin{equation*}\label{eq4.5}
  \tau(J,\chi)=
  -2n_1 \alpha_{K^1}(t)  +\sum_{i=1}^{n_1}\big(  \alpha_{K^1}(\zeta_p^{a_i} t) +
 \alpha_{K^1}(\zeta_p^{-a_i} t)+
\tau(T_{2,q_1},\chi_{a_i})-\tau(T_{2,q_1},\chi_{b_i})
 \big).
 \end{equation*}
    \qed 
\end{prop}

As explained above, given  any metabolizer, one can find a non-trivial character that vanishes  on it.  Therefore, taking $p=q_1$ in Proposition \ref{prop3.2} and applying  the main result of Casson and Gordon~\cite{CG},   one concludes the following.

\begin{corollary}\label{deprel} If a knot $J$ as above is slice and $n_1 > 0$, then for some set of elements $a_i, b_i \in \{0, 1, \ldots , \frac{q_1 -1}{2}\}$, not all 0, the sum 

$$   -2n_1 \alpha_{K^1}(t)  +\sum_{i=1}^{n_1}\big(  \alpha_{K^1}(\zeta_{q_1}^{a_i} t) +
 \alpha_{K^1}(\zeta_{q_1}^{-a_i} t)+
\tau(T_{2,q_{1}},\chi_{a_i})-\tau(T_{2,q_{1}},\chi_{b_i})
 \big) $$ represents $ 0 \in W(\cc(t)) \otimes \zz_{(2)}.$

\end{corollary}

To apply this as an obstruction to knots being slice, we must understand invariants of $W(\cc(t)) \otimes \zz_{(2)}$ better.  This is accomplished in the next section.


\section{Signatures and discriminants}

There are two fundamental types of invariants that can detect the nontriviality of elements  $\tau\in W(\cc(t), \calj) \otimes \zz_{(2)}$: signatures and discriminants. Discriminants can be computed algorithmically (see~\cite{KL1}), but  they take  values in a 2--torsion group, and thus their use in detecting elements of infinite order is quite tricky.  Signatures take value in a torsion free group, $\zz_{(2)}$, but are difficult to compute.  We now describe a method which will allow us to bypass these difficulties by taking advantage of the interplay between signatures and discriminants.  An added advantage of this approach is that it helps us avoid the usually challenging problem of identifying and analyzing all possible metabolizers for the linking forms of the relevant 3--manifolds.   Useful references  for Casson-Gordon discriminant invariants include~\cite{gl2} and~\cite{KL1}.

 \subsection{Basic facts about  $W(\cc(t), \calj) \otimes \zz_{(2)}$.}

 In Appendix~\ref{appwitt} we present some of the details concerning the Witt group $W(\cc(t), \calj) \otimes \zz_{(2)}$.  Here are the key points that we need.

 \begin{itemize}
 
 \item If $I\in W(\cc(t), \calj)$ is represented by a Hermitian matrix $A$  with polynomial entries, the {\em jump function} $$j(I)(\omega)\co S^1 \to \zz$$ represents half the jump in the signature function sign($A(\omega)$) defined for $\omega \in S^1$. The function $j(I)$ has finite support. It extends to a well-defined  $\zz_{(2)}$--valued function on  $W(\cc(t), \calj) \otimes \zz_{(2)}$.
   
 \item The {\em discriminant} of a class $I = [A]  \in W(\cc(t), \calj)$, where the matrix $A$ is of rank $n$,    is given by  
 $$\disc(I)=(-1)^{n(n-1)/2}\det(A).$$ This defines a {\it function} (but not a homomorphism), 
 $$\disc\co W(\cc(t), \calj)\to (\cc(t)^\calj)^\times/N,$$ where $(\cc(t)^\calj)^\times $ denotes the non-zero symmetric   ($f =  \calj(f)$) rational functions and  
 $N $ denotes the {\em norms}; that is, the multiplicative subgroup of    $\cc(t)^\times$    defined as
 $$N= \{f   \calj(f) \ |\  f \in \cc(t), f \ne 0 \}.$$  
 
 \item By taking the further quotient by the subgroup $\pm 1$ there is a well-defined homomorphism $$\disc_\pm \co W(\cc(t), \calj)\to (\cc(t)^\calj)^\times/\pm N.$$
 
 \item A class   $d\in (\cc(t)^\calj)^\times/N$  has a canonical  representative  in $ (\cc(t)^\calj)^\times$   of the form 
 $$d=  at^{-n}\prod_{i=1}^{2n} (t- \omega_i),$$ where the $\omega_i$ are distinct unit complex numbers and $a^2 = 1/ \prod \omega_i$.  If $d = \disc(I)$, then  the set of numbers $\{\omega_i\}$ are called {\it the roots} of $\disc( I)$.
 
 \item There is a natural extension of $\disc_\pm$ to $W(\cc(t), \calj) \otimes \zz_{(2)}$, defined by 
 $\disc_\pm (I\otimes\frac{p}{q}) = (\disc_\pm I)^p$.  This is again a homomorphism.
 
 \item A class $\frac{p}{q}I \in W(\cc(t), \calj) \otimes \zz_{(2)}$ has $j(I)(\omega)$  odd if and only if $ \omega$ is a  root of  $\disc_\pm(\frac{p}{q}I )$.  (An element in $\zz_{(2)}$ is called odd if it is not in $2\zz_{(2)}$ and is called even otherwise.)
 \end{itemize}

  \subsection{Twisted polynomials and the discriminant.}  Let $\chi \co \pi_1(M^3(K)) \to C_p$.  Then, as described in~\cite{KL1}, we may associate to $K$ and $\chi$   a  {\it twisted Alexander polynomial} $\Delta_{K,\chi}(t) \in  \cc[t^\pm]$.  Theorem 6.5 of~\cite{KL1} states:

\begin{theorem} \label{thm:discrim1} 
 $\disc_\pm(\tau(K,\chi))=(1-t)^e\Delta_{{K},\chi}(t),$
where $e=0$ if $\chi$ is trivial and $e=1$ if $\chi$ is non-trivial. \qed
 \end{theorem}
\noindent Note that in this theorem the twisted polynomial is well-defined up to multiplication by $at^k$ for $k \in \zz$ and $a \in \cc^\times$, while the discriminant is well-defined up to $\pm N$. We refer the reader to~\cite[Sections 2.2 and 6]{KL1} for further details, and to~\cite{HKL} for an alternative  description of this twisted Alexander polynomial as a twisted polynomial of a 2--dimensional metabelian representation of $\pi_1(S^3-K)$.
    
Theorem \ref{thm:discrim1} and the discussion that precedes it generalizes to the Casson--Gordon setting   the well-known facts that the discriminant of $\alpha_K(x)$ equals the ordinary Alexander polynomial of $K$ modulo norms, and that the jump function of the Levine--Tristram $\omega$--signatures is supported on the roots of the Alexander polynomial  (see the first paragraphs of Section  \ref{secttorusknot} below).

\subsection{The discriminant and jump function of the torus knot ${\pmb T_{2,p}}$.}
  
As an important example, Theorem~\ref{thm:discrim1} allows us to readily compute the discriminant of the Casson-Gordon invariant of $T_{2,p}$ when $p$ is a prime and $\chi$ is any $C_p$--valued character.  Combined with Corollary~\ref{jumpcor}, we obtain information about the jumps of the  signature function of  $\tau(T_{2,p}, \chi_a)$.

\begin{lemma}\label{discrim} Let $T_{2,p}$  denote the $(2,p)$--torus knot for some odd prime $p$, and $M^3(T_{2,p})$ its 2--fold branched cover.  Let $f(t)=1+t+t^2+\cdots + t^{p-1}$.  

There exists  $d \in\{1,2,\cdots, \tfrac{p-1}{2}\}$ so that for any $a$, 
$$\disc_\pm(\tau(T_{2,p},  \chi_{a}))=    \frac{t^{\frac{3-p}{2}}f(t)}{ (t-\zeta_p^{ad})(t-\zeta_p^{-ad})}. $$
 Hence  if $a\not\equiv 0\pmod{p}$ and $\theta$ denotes the trivial character,
$$  j\big(\tau(T_{2,p},  \chi_{a })-\tau(T_{2,p}, \theta)\big)(\omega)\text{ is }\begin{cases} \text{even}&\text{ if } \omega\ne \zeta^{\pm ad}\\
                                                                                                         \text{odd}& \text{ if } \omega=\zeta^{\pm ad}. \end{cases}$$
\end{lemma}
 
\begin{proof}
    The $(2,p)$--torus knot has presentation $\pi=\langle \alpha, \beta \ | \ \alpha^2\beta^{p}\rangle$.  Define $n=\frac{p-1}{2}$.  The meridian (that is, the generator of $H_1(S^3-T_{2,p})=\zz$) is given by $\mu=\alpha + n\beta$, and in $H_1(S^3-T_{2,p})$, $\alpha=p\mu$ and $\beta=-2\mu$.

We use the methods and notation of~\cite{HKL}.  In that article it is explained how a choice of character $\chi \co H_1(M( T_{2,p}))\to  C_p$ determines and is determined by a dihedral representation of $\pi$. Let $\zz/2=\langle x\ | \ x^2=1\rangle$ act on $C_p=\{\zeta_p^i\}$ via $x\cdot \zeta_p=\zeta_p^{-1}$; then given $d\in \{0,1, \cdots, p-1\}$,  
$$\alpha\mapsto x, \ \beta\mapsto \zeta_p^d$$
determines a $\zz/2\ltimes C_p$  representation since $x^2 = 1 =(\zeta_p^d)^p$.  This representation restricts to a  trivial  representation on the 2--fold cover if and only if $d=0$, since $\beta=-2\mu$ in $H_1(S^3-T_{2,p})$.

From this dihedral representation a recipe is given in~\cite{HKL} to produce a $GL_2(\cc[t^{\pm 1}])$ representation of $\pi$ whose associated twisted Alexander polynomial is $\Delta_{ {K},\chi}(t)$.  The recipe produces the representation $\rho \co \pi\to  GL_2(\cc[t^{\pm 1}])$:
$$\rho(\alpha)= \begin{pmatrix}
0&1\\t&0
\end{pmatrix}^p=t^n\begin{pmatrix}
0&1 \\ t &0\end{pmatrix}, \ \rho(\beta)= \begin{pmatrix} 0&1\\ t&0\end{pmatrix}^{-2}\begin{pmatrix}
\zeta&0\\0&\zeta^{-1}\end{pmatrix}^d=t^{-1}\begin{pmatrix}
\zeta^d&0\\0&\zeta^{-d}\end{pmatrix}.$$

 Theorem 7.1 of~\cite{HKL} shows that $\Delta_{ {K},\chi}(t)$ is   the order of the $\cc[t^{\pm1}]$-torsion of the corresponding twisted first  homology module  $H_1(S^3-K; (\cc[t^{\pm 1}])^2_\rho)$; here $(\cc[t^{\pm 1}])^2_\rho$ is $(\cc[t^{\pm 1}])^2 =  \cc^2 \otimes \zz[t^\pm]$ viewed as a $\zz[\pi_1(M^3_0(T_{2,p}))]$--module via the representation $\rho \otimes \epsilon$, where $\epsilon$ is the canonical action of  $\zz$ on $\zz[t^\pm]$.  Let $\Delta_0$ denote the order of
 $H_0(S^3-K; (\cc[t^{\pm 1}])^2_\rho)$.
 
Note that $H_0(S^3-K; (\cc[t^{\pm 1}])^2_\rho)$  is the cokernel of the matrix obtained by substituting  the extension of $\rho$ to $\zz\pi\to gl_2(\cc[t^{\pm 1}])$   into the matrix
 $$\partial_1= \begin{pmatrix}
 \alpha-1 \\ \beta-1
\end{pmatrix}$$
(this matrix represents the differential on $1$-chains in the universal cover). A simple calculation shows that $H_0(S^3-K; (\cc[t^{\pm 1}])^2_\rho)$ is trivial if $d\ne 0$, and  $\cc[t^{\pm1}]/\langle t-1\rangle$ if $d=0$. Thus
$\Delta_0=(t-1)^{e-1}$ where $e=0$ if $d=0$   and $e=1$ if $d\ne 0$.

 To compute $\Delta_{ {K},\chi}(t)$ we first compute the Fox matrix 
$$ \partial_2=\begin{pmatrix}
1+\alpha & \alpha^2(1+\beta+\cdots+\beta^{p-1})
\end{pmatrix}$$
representing the differential on $2$--chains in the universal cover. Then Theorem 4.1 of~\cite{KL1} shows that 
$$ \Delta_{ {K},\chi}(t) =\frac{\det(\rho(1+\alpha))}{\det(\rho(\beta-1))}\Delta_0=\frac{\det(\rho(  \alpha^2(1+\beta+\cdots+\beta^{p-1})))}{\det(\rho(\alpha-1))}\Delta_0.$$

Now 
$$\det(\rho(1+\alpha))=\det \begin{pmatrix}
1&t^n  \\
t^{n+1} &1 
\end{pmatrix}=1-t^p
$$
and
$$\det(\rho(\beta-1))=\det \begin{pmatrix}
t^{-1}\zeta_p^d-1&0  \\
0&t^{-1}\zeta_p^{-d}-1
\end{pmatrix}= t^{-2}(t-\zeta_p^d)(t-\zeta_p^{-d}). 
$$
 Using Theorem~\ref{thm:discrim1} we find that for some $a$ and $k$, 
 $$\disc_\pm(\tau(T_{2,p},\chi))=a t^k (1-t)^e\Delta_{ {K},\chi}(t)=a t^k (1-t)^{2e-1}t^2\frac{1-t^p}{(t-\zeta_p^d)(t-\zeta_p^{-d})}.
$$
 Since $-t^{-1}(1-t)^2=(1-t)\calj(1-t)$,  and $f(t)(t-1)=t^p-1$, this can be rewritten (perhaps changing $a$ and $k$) as 
$$\disc_\pm(\tau(T_{2,p},\chi))= a t^k \frac{f(t)}{(t-\zeta_p^d)(t-\zeta_p^{-d})}.
$$
 Symmetry of the discriminant then implies that $k =  \frac{3-p}{2}$ and $a = \pm 1$.

The lemma follows from the fact that all non-trivial characters are multiples of $\chi_1$. Hence if $d$ is chosen so that the character $\chi_1$   takes $\beta$, viewed as a loop in the 2--fold cover, to $\zeta_p^d$, $\chi_a$ corresponds to the character that takes $\beta$ to $\zeta_p^{ad}$. 

\end{proof}

Given an odd prime $p$, define a homomorphism 
\begin{equation}
\Psi_p \co W(\cc(t), \calj)\otimes \zz_{(2)}  \to ( \zz_{(2)})^{(p-1)/2}
\end{equation}
by evaluating the jump function at the non-trivial $p$--roots of 1 in the upper half-circle: 
$$\Psi_p(I)=  \big(j(I)(\zeta_p), j(I)(\zeta_p^2),\cdots, j(I)(\zeta_p^{(p-1)/2})\big).$$
Note that      for $I\in W(\cc(t))\otimes \zz_{(2)}$, 
 $  \Psi_p(I)\in (2\zz_{(2)})^{(p-1)/2}$ if and only if $\disc_\pm(I)$ has no roots among   $\zeta_p,\zeta_p^2,\cdots, \zeta_p^{(p-1)/2}$.

\begin{corollary}\label{cor5.5} If $p$ is an odd prime, the set of Witt classes
$$ \tau(T_{2,p},  \chi_{a })-\tau(T_{2,p}, \theta)\in W(\cc(t), \calj)\otimes   \zz_{(2)}, \ a=1,2,\cdots,\tfrac{p-1}{2}$$
are linearly independent  and their span is mapped injectively to $  (\zz_{(2)})^{(p-1)/2}$ by $\Psi_p$.
\end{corollary}
\begin{proof}   Consider the homomorphism $$\Phi \co (\zz_{(2)})^{(p-1)/2}  \to W(\cc(t), \calj)\otimes  \zz_{(2)}$$ which takes  the $a$th coordinate vector to the difference  $\tau(T_{2,p}, \chi_a)-\tau(T_{2,p},\theta)$ in $ W(\cc(t), \calj)\otimes  \zz_{(2)}$.

Lemma \ref{discrim}  implies that the matrix for $p \cdot \Psi_p \circ \Phi$ differs from a permutation of the identity by  an even matrix, and hence has odd (and, in particular, non-zero) determinant.  It follows that $\Phi$ is injective. \end{proof}

 \section{The main examples}
 
 In Section~\ref{sectionli} we considered the knots $$J_i=  K^i _{2,q_{2i-1}} + T_{2,  q_{2i}} -   K^i_{2, q_{2i}} - T_{ 2, q_{2i-1}}.$$   Our goal is to prove that for appropriate choices of $K^i$ and $q_i$, the set $\{J_i\}_{i=1}^\infty$ is linearly independent.

 The conditions on the knots $K^i$   which we will need to arrive at a contradiction (to the assumption that $J$ is slice) are that  the $K^i$ be  $p$--{\it deficient} and $p$--{\it independent}.   These are conditions on the algebraic concordance class of $K^i$.  Roughly stated, $K$ is $p$--deficient if its (Levine--Tristram)  signature function has no jumps at $p$th roots of unity and is $p$--independent if the abelian Witt invariant $\alpha_K(t)$ and its translates $\alpha_K(\zeta_p^at)$ are linearly independent in $W(\cc(t))$.  
 
 \begin{definition}
   Given a knot $K$ and an odd prime $p$,  we say that $K$ is {\it $p$--deficient}  if     $j(\alpha_K(t))(\zeta_p^a) = 0$ for all $a\in\{0,1,\cdots, p-1\}$. 
\end{definition}
 
  \begin{definition} Given a knot $K$ and an odd prime $p$,  we say that $K$ is {\it $p$--independent} if  the  Witt classes $ \alpha_K(\zeta_p^a t), $ $a\in\{0,1,\cdots, p-1\}$, in $W(\cc(t)) \otimes \zz_{(2)}$ are linearly independent.
 \end{definition}

\begin{lemma}\label{defindeplem} If a knot $K$ is $p$--deficient and $p$--independent, then for any choice of integers $n>0$ and $a_1,\cdots, a_n\in\{0,1,\cdots, p-1\}$ with not all $a_i$ zero, 
$$   -2n  \alpha_{K }(t)  +\sum_{i=1}^{n }\big(  \alpha_{K }(\zeta_{p}^{a_i} t) +
 \alpha_{K }(\zeta_{p}^{-a_i} t)\big)$$
 is a non-zero element of the kernel of $\Psi_p\co W(\cc(t)) \otimes \zz_{(2)}\to (\zz_{(2)})^{(p-1)/2}$.
\end{lemma}

\begin{proof} Note that $j(\alpha_K(\zeta_p^{a_i}))(\zeta_p^a)=j(\alpha_K)(\zeta_p^{a_i+a})$, which vanishes since $K$ is $p$--deficient.  Hence $\Psi_p\big(    -2n  \alpha_{K }(t)  +\sum_{i=1}^{n }\big(  \alpha_{K }(\zeta_{p}^{a_i} t) +
 \alpha_{K }(\zeta_{p}^{-a_i} t)\big) \big)=0$.
 
 Since  $K$ is $p$--independent and  $-2n  \alpha_{K }(t)  +\sum_{i=1}^{n }\big(  \alpha_{K }(\zeta_{p}^{a_i} t) +
 \alpha_{K }(\zeta_{p}^{-a_i} t)\big)$ is a non-trivial  (not all $a_i$ are zero) linear combination of the  $\alpha_K(\zeta^at)$, it is non-zero. 
\end{proof}

 Nontrivial examples of $p$--deficient and $p$--independent   knots will be presented in Section~\ref{secttorusknot}. In particular, we will show that the trefoil, $T_{2,3}$, is $p$--deficient and $p$--independent  for all primes $p>3$.

 \begin{theorem} Let $J=\sum_{i=1}^N n_iJ_i$ with the $J_i$ as above.  If, for some $j$ with $n_j\ne 0$,  the knot $K^j$ is $q_j$--deficient and $q_j$--independent, then $J$ is not slice.
\end{theorem}

\begin{proof}
  Suppose that $J$ is slice.  Assume, by changing sign and reindexing if necessary, that  $j=1$ and that $n_1>0$.

 In this case, we found in Corollary~\ref{deprel} that 
    for some set of elements $a_i, b_i \in \{0, 1, \ldots , \frac{q_1 -1}{2}\}$, not all 0, then 
$$   -2n_1 \alpha_{K^1}(t)  +\sum_{i=1}^{n_1}\big(  \alpha_{K^1}(\zeta_{q_1}^{a_i} t) +
 \alpha_{K^1}(\zeta_{q_1}^{-a_i} t)+
\tau(T_{2,q_{2i-1}},\chi_{a_i})-\tau(T_{2,q_{2i-1}},\chi_{b_i})
 \big) =0$$ in $  W(\cc(t)) \otimes \zz_{(2)}.$

Applying the function $\Psi_{q_1}$ to this equation we find, using Lemma \ref{defindeplem}, that $$\Psi_{q_1}(\tau(J,\chi))=\Psi_{q_1}\big(\sum_{i=1}^{n_1}\tau(T_{2,q_1},\chi_{a_i})-\tau(T_{2,q_1},\chi_{b_i})
 \big) = 0.$$ 
  This   can be rewritten as 
  $$\Psi_{q_1}\big(\sum_{i=1}^{n_1}\big(\tau(T_{2,q_1},\chi_{a_i}) - \tau(T_{2,q_1},\theta)\big) -\big(\tau(T_{2,q_1},\chi_{b_i}) - \tau(T_{2,q_1},\theta) \big)
 \big) = 0.$$ 
 By Corollary~\ref{cor5.5}, this implies that   $$ \sum_{i=1}^{n_1}\big( \tau(T_{2,q_1},\chi_{a_i}) - \tau(T_{2,q_1},\theta)\big) -\big(\tau(T_{2,q_1},\chi_{b_i}) - \tau(T_{2,q_1},\theta) \big)
 = 0,$$ and thus 
 $\sum_{i=1}^{n_1}\tau(T_{2, q_1},\chi_{a_i})-\tau(T_{2,q_1},\chi_{b_i})
=0$.  We  also conclude that the  (unordered) sets $\{a_1,a_2,\cdots, a_{n_1}\}$ and $\{b_1,b_2,\cdots, b_{n_1}\}$ coincide.   In particular, at least one of the $a_i$ is non-zero.  Thus
\begin{equation*}
\label{eq6.2}
0=-2n_1 \alpha_{K^1}(t)  +\sum_{i=1}^{n_1}  \alpha_{K^1}(\zeta^{a_i } t) +
 \alpha_{A_1}(\zeta^{-a_i} t). 
\end{equation*}
But this is impossible by Lemma \ref{defindeplem}.
Hence $J$ cannot be slice. 
 \end{proof}
 
 With this, our main result follows.  
 
 \begin{theorem}\label{maintheorem}  Let $q_i$ be a sequence of  positive integers with $q_{2i - 1}$ prime for all $i$ 
 and $q_{2i}$ relatively prime to $q_{2j-1}$ for all $i,j$.
Let $K^i$ be a sequence of knots so that $K^i$ is $q_{2i-1}$--deficient and $q_{2i-1}$--independent for all $i$.   Let $$J_i = K^i_{ 2, q_{2i-1}} \cs T_{ 2, q_{2i}} \cs -K^i_{2,  q_{2i}} \cs -T_{ 2, q_{2i-1}}.$$  Then the $J_i$ are linearly independent, algebraically slice knots.  \qed

 \end{theorem} 
 
 As mentioned above, the next section shows that $T_{2,3}$ is both  $p$--deficient and $p$--independent  for all primes $p>3$.  Given this, the following corollaries are immediate, and yield Theorems 1, 2, and 3 of the introduction.
 
 \begin{corollary}\label{maincor}
The algebraically slice  knot $$T_{2,3;2, 13} + T_{2,15}-T_{2,3;2,15} -T_{2,13}$$ has infinite order in $\calc$.  
\end{corollary}
\begin{proof}  The assertion   follows immediately from Theorem \ref{maintheorem}.
\end{proof}

\begin{corollary}\label{large} Let $q_1=13, q_2=17, \cdots$  be the increasing list of primes   greater than 11. Then the set of algebraic knots 
 $$\{ T_{2,q_i}, \ T_{2,3;2,q_i}\}_{i=1}^\infty$$
 is a basis for a free abelian subgroup of  the concordance group $\calc$. This subgroup intersects the kernel of $\calc\to\calg$ in a free abelian subgroup, with basis the set of algebraically slice knots
 $$ \{T_{2,q_i} -T_{2,3;2,q_i}-T_{2,13} +T_{2,3;2,13}\}_{i=2}^\infty.$$
\end{corollary}
\begin{proof}

Consider a linear combination
$$J=\sum_{i=1}^N n_i T_{2,q_i} + m_i T_{2,3;2,q_i}.$$
Suppose that $J$ is slice. We will show that each $n_i$ and $m_i$ is zero.

Fix $\ell$.  When evaluated at $\omega=e^{2\pi i/(2q_\ell)}$, 
the jump function for the Levine--Tristram signature of a knot in $\{ T_{2,q_i}, \ T_{2,3;2,q_i}\}_{i=1}^\infty$ is non-zero only for the knots $T_{2,q_\ell}$ and $T_{2,3;2,q_\ell}$, since the $q_i$ are different primes.  Indeed, for $T_{2,q_\ell}$ and $T_{2,3;2,q_\ell}$, the jump is equal to $-1$  (see, for example,~\cite{Litherland1979}).  This implies that $m_\ell=-n_\ell$.

Furthermore, the  jump function for the Levine-Tristam signature of $T_{2,3;2,q_i}$, evaluated at $\omega=e^{2\pi i/12}$, is equal to $-1$ for  all $i$.   For this value of $\omega$, however, the jump function for  $T_{2,q_i}$ is zero.    It follows that the sum of the $n_i$ is zero.

Thus \begin{equation}\label{eq6.1}
J=\sum_{i=1}^N n_i  (T_{2,q_i} -T_{2,3;2,q_i}) \text{ with } \sum_{i=1}^N n_i=0. 
\end{equation}
Any knot of the form (\ref{eq6.1}) is algebraically slice: indeed, its Blanchfield form is 
$$Bl_J(t)= \sum_{i=1}^N n_i( Bl_{T_{2,q_i}}(t) - Bl_{T_{2,3}}(t^2)-  Bl_{T_{2,q_i}}(t))=-\sum_{i=1}^N n_i(  Bl_{T_{2,3}}(t^2))=0.$$

Since $\sum n_i=0$, we have
$$J=\sum_{i=1}^N n_i  (T_{2,q_i} -T_{2,3;2,q_i})- \sum_{i=1}^N n_i  (T_{2,11} -T_{2,3;2,11}),$$
as an equation in $\calc$. Theorem \ref{maintheorem}, together with the fact that $T_{2,3}$ is $p$--deficient and $p$--independent for $p>3 $, implies that each $n_i$ is zero.  This proves that the set $\{T_{2,q_i}, T_{2,3;2,q_i}\}$ is linearly independent.

Since the   jumps in the Levine--Tristram signature functions are determined by the algebraic concordance class of a knot, (\ref{eq6.1})  shows that  the intersection of the span of  $\{T_{2,q_i}, T_{2,3;2,q_i}\}$ with the kernel of $\calc\to \calg$ is  a free abelian group,  with basis the set  of algebraically slice knots $\{T_{2,q_i} -T_{2,3;2,q_i}-T_{2,13} +T_{2,3;2,13}\}_{i=2}^\infty$.

\end{proof}

\begin{corollary}
Let $\cala \subset \calc$ denote the subgroup of the knot concordance group generated by algebraic knots. 
The intersection of $\cala$ with the kernel of the map $\calc \to \calg$ to the algebraic concordance group contains an  infinitely generated free abelian group.\qed
\end{corollary}

\section{Torus knot examples:  $p$--deficiency and $p$--independence.}\label{secttorusknot}

Let $K$ be a knot in $S^3$, $F$ a Seifert surface for $K$ and $V$ the Seifert form for $F$.  There are several well-known constructions of  a $4$-manifold $X^4$ with boundary $S_0^3(K)$ over which $\delta\co H_1(S^3_0(K))\to \zz$ extends such that the equivariant intersection form of $X^4$ is $I(X^4,\bar{\delta})=(1-x)V+(1-x^{-1})V^T$ and   intersection form $I(X)=(1)$.  Such constructions can be found in~\cite{CG, Ka, KKR}.

It follows that $\alpha_K(x)\in W(\qq(t))$ is represented   by the matrix
$$\begin{pmatrix}(1-x)V+(1-x^{-1})V^T&0\\0&-1
\end{pmatrix}.$$

Since $(x^{-1}-1)\big(xV-  V^T\big)=(1-x)V+(1-x^{-1})V^T$, and the Alexander polynomial satisfies
$$\Delta_K(x)=\det\big(xV-  V^T\big),$$
it follows that the jumps in the Levine-Tristram signature function of $\alpha_K(x)$ is supported on a subset of the roots 
 of the Alexander polynomial.  Notice that this is a more precise statement than saying that the odd jumps occur at roots of the discriminant, since the Alexander polynomial is well--defined in  $\zz[x^{\pm1}]$. Furthermore, if $\omega$ is a root of unity and  $(x-\omega)$ divides $\Delta_K(x)$ but $(x-\omega)^2$ does not divide $\Delta_K(x)$, then 
 $j_\omega(\alpha_K)=\pm 1, $  improving the conclusion of Corollary \ref{corA6}.

\begin{theorem}For any relatively prime integers $m$, $n$, and $q$, and any prime divisor $p$ of $q$, the torus knot $T_{m,n}$ is $p$--deficient and $p$--independent.

\end{theorem}

\begin{proof}
The Alexander polynomial of $T_{m,n}$ is $$\Delta_{T_{m,n}}(x) = \frac{(x^{mn} -1)(x-1)}{(x^m-1)(x^n-1)}.$$ Thus, the only roots of $\Delta_{T_{m,n}}(x)$ are the $mn$--roots of unity which are not simultaneously $m$ or $n$--roots of unity.  It follows that the jumps in the Levine--Tristram signature function of $\alpha_{T_{m,n}}$ occur, and equal $\pm1$, at the unit complex numbers $e^{2\pi i\frac{c}{mn}} = \zeta_{mn}^c$, where $c$   is not divisible by either $m$ or $n$, and $1 \le c \le mn-1$.  (There are $(m-1)(n-1)$ such $c$.)
\vskip.1in
\noindent{\bf $p$--deficiency:}   From the definition, we see that if $T_{m,n}$ is not $p$--deficient, then for some $a\in \{0,1,\cdots, p-1\}$ and $c$ as in the previous paragraph,   $\zeta_p^a = \zeta_{mn}^c $. This is impossible, since $p$ and $mn$ are relatively prime and $1\leq c \le mn -1$.  

\vskip.1in

\noindent{\bf $p$--independence:} To demonstrate  the independence of the $\alpha_{T_{m,n}}(\zeta_p^at)$,  we show that for distinct $a_1$ and $a_2$,  $0 \le a_1, a_2 \le p-1$, the jumps for the Levine--Tristram signature function occur at distinct points. The jumps for $\alpha_{T_{m,n}}(\zeta_p^{a_i}t)$ occur at $\omega=\zeta_p^{-a_i}\zeta_{mn}^{c_i}$, where $c_i$ is not divisible by either $m$ or $n$ and $1 \le c_i \le mn-1$.

If the jumps occured at the same point, then $\zeta_p^{-a_1}\zeta_{mn}^{c_1}=\zeta_p^{-a_2}\zeta_{mn}^{c_2}$, and so 
$$ \frac{c_1}{mn} - \frac{a_1}{p} = \frac{c_2}{mn} - \frac{a_2}{p} \mod \zz.$$  This can be rewritten as:
 $$\frac{ (c_1-c_2)p  - (a_1 - a_2)mn}{mnp} \in \zz.$$
 This immediately implies that $a_1 - a_2$ is divisible by $p$, which in turn implies that $a_1 = a_2$, giving the desired contradiction.

\end{proof}

 
\section{The 4--ball genus.}

 We next observe that  if $q_1,q_2$ are a pair of integers and $K$ is any knot, the algebraically slice knots 
$$ J =  K_{2,q_1}  -   K_{2,q_2} - T_{2, q_1} +T_{2, q_2}$$  have 4--ball genus equal to $0$ or $1$.  In the case that $K$ is slice, we noted in Lemma~\ref{summarylemma} that $J$ is slice.  If 
 $K$ is $q_1$--deficient and independent, then Corollary~\ref{maincor} shows $J$   is not slice.  The following argument shows that in this second case   $J$ has 4--ball genus at most 1.

\begin{figure}
\psfrag{A}{\huge{$K$}}
\psfrag{-A}{\huge{$-K$}}
\psfrag{k1}{$q_1$}
\psfrag{k2}{$q_2$}
\psfrag{-k1}{$-q_1$}
\psfrag{-k2}{$-q_2$}
\psfrag{a1}{$\gamma_1$}
\psfrag{a2}{$\gamma_2$}
\psfrag{a3}{$\gamma_3$}

\begin{center}
\includegraphics[height=200pt]{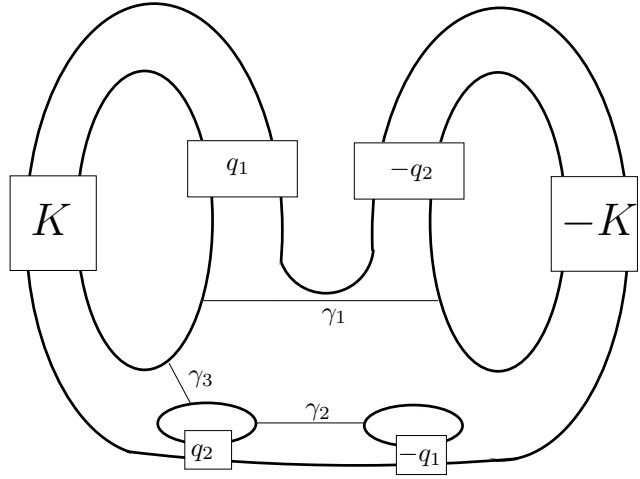}
\caption{\label{fig:3}$ J =  K_{2,q_1}  -   K_{2,q_2} - T_{2, q_1} +T_{2, q_2}$}
\end{center}
\end{figure}

Figure~\ref{fig:3} illustrates $J$   with three arcs,  $\gamma_1,\gamma_2,$ and $\gamma_3$, depicted. In this figure the labels $\pm q_i$ refers to half-twists, and the parallel strands  passing through $\pm K$ are to be tied in the knot $\pm K$.

 A band move along $\gamma_1$ produces a satellite (link) of the slice knot $K\#(-K)$.   Taking the corresponding satellite of the null-concordance and then performing band moves along the arcs labelled $\gamma_2$ and $\gamma_3$ gives a genus 1 cobordism from $J(K,k_1,k_2)$ to a 2-component unlink.   Thus $J(K,k_1,k_2)$ has 4--ball genus at most one.

\begin{corollary} The  knot  $T_{2,3;2, 13} + T_{2,15}-T_{2,3;2,15} -T_{2,13}$ has 4-ball genus equal to 1.\qed
\end{corollary}

\vskip.1in

As mentioned in the introduction, the Ozsv{\'a}th-Szab{\'o} and Rasmussen concordance invariants, $\tau$ and $s$, are unable to determine whether any of the algebraically slice linear combination involving positive iterated torus knots, is slice.  We make this precise in the following proposition.

\begin{prop}\label{prop:stau} Fix $p,q_1,q_2>0$.  Suppose $J =  K_{p,q_1}  -   K_{p,q_2} - T_{p, q_1} +T_{p, q_2}$, with $K=K_{r_1,s_1;\cdots;r_n,s_n}$ a positively iterated torus knot; that is, $r_i,s_i>0$ for all $i$.  Then $\tau(J)=s(J)=0,$ where $\tau$ and $s$ are the  Ozsv{\'a}th-Szab{\'o} and Rasmussen concordance invariants, respectively.
\end{prop}
\begin{proof}  It is well-known that the Seifert genus of $K_{p,q}$ is given by $$g(K_{p,q})=pg(K)+ g(T_{p,q}),$$ and that $2g(T_{p,q})=(p-1)(q-1)$. See, for instance,~\cite[Chapter 1\S 3]{EN1985}. 

We claim $2\tau(K)=s(K)=2g(K)$ for any positively iterated torus knot.  Given this, the proposition follows from the genus formula above, together with the fact  that both invariants change sign under reflection, $\tau(K)=-\tau(-K)$ and $s(K)=-s(-K)$.

For torus knots, the fact that $2\tau(K)=s(K)=2g(K)$ was proved by Ozsv{\'a}th and Szab{\'o}~\cite{Lens} and Rasmussen~\cite{ra}, respectively.  For positively iterated torus knots, the result follows 
from~\cite[Corollary 1.4]{ComplexCable}, which  shows that a positively iterated torus knot bounds a Seifert surface which is isotopic to a piece of a complex curve in the four-ball (for algebraic knots this is well-known, through the work of Milnor~\cite{Milnor1968}).   Knots which bound such complex curves satisfy the stated equalities, by~\cite{Livingston2004} (see also~\cite{SQPfiber}).
\end{proof}


\appendix


\section{Properties of the Witt group $W(\cc(t))$}\label{appwitt}

\subsection{Hermitian forms over $\cc  (t) $ and the Witt group $W(\cc(t), \calj)$.}  The ring $\cc[t^\pm]$ has the involution denoted $\calj$, defined by  $\calj({\sum a_i t^i}) = \sum \overline{a}_i t^{-i}$, where $\overline{a}_i$ denotes standard complex conjugation.  There is a natural extension of $\calj$ to $\cc(t)$ and then to $GL_n(\cc(t))$, applying $\calj$ to the entries of a matrix and transposing.   

A Hermitian inner product space consists of a pair $(V, \beta)$ where $V$ is a finite dimensional vector space over $\cc(t)$ and $\beta$ is a Hermitian inner product:  $\beta(av, bw) = a \calj({b})(v,w)$ and $\beta(v,w) = \calj({\beta(w,v)})$.  In terms of a basis, the Hermitian inner product  $\beta$ is given by a matrix $B$ satisfying $\calj (B)=B$. Given a second basis, the new coordinates are related to the old by a change of basis matrix $P$, and the new matrix representation of $\beta$ is $PB\calj({P})$. 
 The inner product is called non--singular if $\det(B) \ne 0$ for some (and thus any) matrix representation $B$ of $\beta$.

An inner product space is called {\it Witt trivial} if there is a half-dimensional subspace of $W \subset V$ such that
$\beta(w_1, w_2) = 0$ for all $w_1, w_2 \in W$.  Two nonsingular Hermitian inner product spaces $(V_1, \beta_1)$ and $(V_2, \beta_2)$ are called {\it Witt equivalent} if $(V_1 , \beta_1) \oplus (V_2, -\beta_2)$ is Witt trivial.  The set of equivalence classes forms an abelian group, with addition induced by direct sum.  This group is the Witt group, which we denote $W(\cc(t))$.  Note that isometric forms are Witt equivalent.

\subsection{Signature invariants defined on $W(\cc(t),\calj)$.}  Let $A(t)$ be a nonsingular Hermitian matrix over $\cc(t)$ representing a class in $W(\cc(t), \calj)$.  For any unit complex number $\omega$, the complex matrix $A(\omega)$ is Hermitian.  The matrix may be singular for isolated values of $\omega$, but in any case it has a signature,  $s_ \omega(A) \in \zz$.  This does not yield a  well-defined homomorphism  $W(\cc(t),\calj)\to \zz$, because of the possible singularities; the problem is corrected by averaging the one-sided limits, defining 
$$\sigma_{e^{i\alpha}}([A])=\tfrac{1}{2}\big( \lim_{x\to\alpha+}  s_{e^{ix}}(A) +  \lim_{x\to\alpha-}  s_{e^{ix}}(A) \big).$$
This results in a homomorphism 
$$\sigma \co W(\cc(t), \calj)\to \text{Funct}(S^1,\zz),\ \ \ \sigma(I)(\omega)=\sigma_{\omega}(I),$$ 
which extends   in the obvious way to 
 $$\sigma \co W(\cc(t))\otimes R\to \text{Funct}(S^1,R)$$
 for any subring $R\subset \qq$.  In our applications we will take $R= \zz_{(2)}$, the set of rational numbers with odd denominator. 

 It is most convenient to re-express this invariant  in terms of the {\em signature jump function}
$$j  \co  W(\cc(t)) \to \text{Funct}(S^1,\zz),\ \ j(I)(e^{i\theta})=
\tfrac{1}{2}\big(\lim_{x\to\theta+}  \sigma(I)(e^{ix})-\lim_{x\to\theta -}  \sigma(I)(e^{ix})\big).$$
We have divided by two to avoid confusion below. This jump function has finite support and  is integer valued, as follows from a diagonalization argument and an examination of 1--dimensional forms.  Details appear later in this appendix.

Since $j(k I)=k j (I)$ for any integer $k$, the signature jump function extends to a homomorphism
$$j \co  W(\cc(t))\otimes R\to \text{Funct}(S^1, R),$$
for any subring $R$ of $\qq$.   

\subsection{Discriminants.}  Roughly stated, the discriminant of a  class in $W(\cc(t), \calj)$ is given by the determinant of a matrix representative of the class.  Since the determinant of a Witt trivial matrix is of the form $\pm f(t)\calj(f(t))$, one needs to view it in a quotient of $\cc(t)$ by such norms.  Taking care in the treatment of signs leads to the formal definition.  For $A$, a matrix of dimension $k$ representing a class $[A] \in W(\cc(t),\calj)$, define
 $$\disc([A]) = (-1)^{k(k-1)/2} \det(A) \in (\cc(t)^\calj)^\times/ N(\cc(t)),$$ where $(\cc(t)^\calj)^\times$ is the multiplicative  subgroup of  $\cc(t)^\times$ consisting of symmetric
 non-zero rational functions   $$(\cc(t)^\calj)^\times= \{ f(t) \ne 0  \in \cc(t) \ |\  \calj(f(t)) = f(t), f(t) \ne 0\} ,$$  and   $N(\cc(t)) $ is the subgroup of {\em norms}
 $$N(\cc(t))= \{ f(t)\calj(f(t))| f\in \cc(t) \} .$$

Since $\frac{1}{f(t)}$ is equivalent to $\calj(f(t))$ modulo norms and $\cc$ is algebraically closed, 
any class in $(\cc(t)^\calj)^\times/ N(\cc(t))$ can be expressed as a factored Laurent polynomial 
\begin{equation}\label{eqapend}
F(t) = at^k \prod (t- \omega_i),
\end{equation} as we now show.

  Since $\calj(F(t)) = F(t)$, it follows that if $(t-\omega)$ is a factor, then $(t - (\bar{\omega})^{-1})$ is also a factor. 
 If $\omega$ is not a unit complex number, then these two factors are distinct.  As both factors appear and $$(t-\omega)(t - (\bar{\omega})^{-1})=(-  \bar{\omega} ^{-1}t) (t-\omega)\calj(t-\omega),$$
they can be removed from the product, modulo norms, at the cost of changing $a$ and the exponent of $t$ in (\ref{eqapend}). 

If $\omega$ is a unit complex number, then 
$$(t-\omega)\calj(t-\omega)=-t^{-1}\bar{\omega}(t-\omega)^2.$$

 Thus we can also remove factors of the form $(t-\omega)^2$  in (\ref{eqapend}) at the cost of changing $a$ and the exponent of $t$.   Any positive real number is a norm, so we may assume $a$ is a unit complex number. A final application of symmetry restricts the form further, yielding the following:

\begin{theorem} Every element in  $(\cc(t)^\calj)^\times/ N(\cc(t))$ has a canonical representative of the form $$at^{-n }\prod_{i=1}^{2n}(t- \omega_i)$$ where the $\omega_i$ are distinct unit complex numbers and $a^2 = 1/ \prod \omega_i$.\qed
\end{theorem}

\begin{corollary}\label{exponent} Let $I \in  W(\cc(t),\calj)$.  For each unit complex number $\omega $, the exponent of $(t-\omega) $ in $\disc(I)$ gives a well-defined homomorphism   $\delta_\omega \co W(\cc(t)),\calj) \to \zz/2\zz$. This extends to a $\zz/2\zz$--valued homomorphism on $W(\cc(t),\calj) \otimes \zz_{(2)}$. \qed
\end{corollary}

\vskip.1in
\subsection{\bf Relation between the discriminant and jump function. }  Via diagonalization, any class in $W(\cc(t), \calj)$ can be represented by a diagonal matrix with diagonal entries in $(\cc(t)^\calj)^\times$.  Scaling basis  elements changes a diagonal entry by an arbitrary norm, and so the preceding discussion applies equally to the diagonal entries.  Thus:

\begin{theorem}\label{diagthm} Every nonsingular Hermitian inner product on $\cc(t)$ has a diagonal matrix representation with each diagonal entry of the form $d = at^{-n} \prod_{i = 1}^{ 2n}(t-\omega_i)$ for some set of distinct  $\omega_i \in S^1$, where $a^2 = 1/\prod \omega_i $.  \qed
  \end{theorem}

A simple trigonometric calculation yields:

\begin{lemma}  If $d =  at^{-n} \prod_{j= 1}^{2n}(t-\omega_j) $, $\omega_j=e^{i\theta_j}$  and $t = e^{i\theta}$, then $$d=(-4)^n \prod_{j=1}^{2n} \sin ((\theta - \theta_j)/2 ).$$\qed
\end{lemma}

This allows one to compute the jump function near $\omega_i$, since $\sin((\theta-\theta_i)/2)$ is negative for $\theta<\theta_i$ and positive for $\theta>\theta_i$.  Together with Corollary \ref{exponent} this implies:

\begin{lemma} \label{lemA5} Let $I \in W(\cc(t), \calj)$ be a 1--dimensional form with canonical representative $ \big(at^{-n} \prod_{i = 1}^{2n}(t-\omega_i) \big)$  as above.   Then 
\begin{enumerate}
\item $\sigma_\omega(I) =   0 $ or $  \pm 1$ depending on whether  $\omega =  \omega_i$ for some $i$ or not.
 \item $j_\omega(I) = \pm 1$ or $0$  depending on whether  $\omega = \omega_i$ for some $i$ or not.

\item $\delta_\omega(I) = 1$ or $0 \in \zz/ 2\zz$, depending on whether $\omega = \omega_i$ for some $i$ or not.
\end{enumerate}\qed
\end{lemma}
 
Since each of these functions is additive and all forms can be diagonalized, we have the following corollary.

\begin{corollary} \label{corA6} 
  Let $ I \in W(\cc(t), \calj)$ be a class with discriminant $d = at^k\prod (t - \omega_i)$, with all   $\omega_i$ distinct unit complex numbers.  Then the jump $j_\omega(I)$ is $0$ or $1 \mod 2$  depending on whether $t-\omega$ is a factor in the discriminant.\qed
 
\end{corollary}

\subsection{Tensoring with $\zz_{(2)}$.}  If $R$ is any subring of $\qq$, then the signature and jump function clearly extend to the Witt group $W(\cc(t),\calj)\otimes R$.  This is not true for the discriminant.  However,

\begin{theorem} After taking the further quotient with $\pm1$, the discriminant extends to a homomorphism $\disc_\pm \co W(\cc(t),\calj)\otimes \zz_{(2)} \to (\cc(t)^\calj)^\times/ \pm N(\cc(t))$.
\end{theorem}
\begin{proof}Any class in $W(\cc(t),\calj)\otimes \zz_{(2)}$ can be expressed as $I \otimes \frac{p}{q}$ with $q$ odd, $p \ge 0$, and $\gcd(p,q) = 1$. We define the discriminant of such a form to be $\disc_\pm (p I) = (\disc_\pm(I))^p$.  Notice however that since $\disc_\pm$ takes values in a group that consists of 2--torsion, the assumption that $\gcd(p,q)=1$ is not needed; introducing odd factors in the numerator does not affect the value of  $(\disc_\pm(I))^p$.  To check this is well defined on the tensor product, one needs to show its value is the same on $kI \otimes \frac{p}{q}$ and $ I \otimes \frac{kp}{q}$ for any integer $k$.  For $k$ even, both are $1$ modulo $\pm N$.  For $k$ odd, both equal $\disc(I)$ modulo $\pm N$.  \end{proof}  

\begin{corollary}\label{jumpcor}
  Let $ \frac{p}{q}I \in W(\cc(t), \calj) \otimes \zz_{(2)}$.  Then the jump $j(I)(\omega)$ is even or odd  (that is, equals $0$ or $1  \in   \zz_{(2)}/ 2 \zz_{(2)}$) depending on whether or not $ \omega$ is a root of the discriminant in its canonical form. \qed
 
\end{corollary}


\section{Ribbon-slice}

\subsection{Connections to the geometry of the ribbon-slice problem}
While the problem of concordance relations between algebraic knots grew out of, and has obvious interest within,  the classical theory of plane curve singularities, it could potentially be more important to the old question of whether slice knots are ribbon.  Indeed, these connections were the original motivation for the work at hand, as we now briefly explain.

First, recall that a knot is {\it ribbon} if it bounds a  {\it ribbon disk} in the four-dimensional unit ball.  A ribbon disk is a smoothly, properly, embedded disk on which the radial function for the four-ball restricts to a Morse function with no index $2$ critical points; that is,  a disk for which the radial function has no local maxima.  The ribbon-slice problem asks whether all slice knots bound a ribbon disk (see for instance~\cite{k}). Apart from the fact that it is longstanding, this question is important because the property of being ribbon has an entirely three-dimensional interpretation. (The interpretation is in terms of immersed disks in the three-sphere with certain types of double points, of so-called ribbon type.)  An affirmative answer would  indicate that the elusive four-dimensional nature of the smooth concordance group could be understood by three-dimensional techniques. 

An interesting but elementary observation  which first appeared in~\cite{Hass}  is that a geometrically minimal disk, that is,  a disk whose mean curvature vector vanishes identically, is a ribbon disk.  This follows from the well-known fact that the  coordinate functions on $\rr^4$ restrict to harmonic functions on a minimal disk. From this it  follows easily that the radial function restricts to a sub-harmonic function, and hence it has no local maxima.
 
In light of this, one is naturally led to wonder  whether  slice knots bound minimal disks in the four-ball.  This is the realm of the classical Plateau problem which asks, more and less generally, when a given curve in $\rr^n$ arises as the boundary of an area-minimizing map  $u \co D^2\rightarrow \rr^n$.  In this form, the answer is ``always," as proved by  Douglas and  Rad{\'o} independently in the early 1930's~\cite{Douglas,Rado}.  The astute reader should be confused, since this seems to imply that not only is every slice knot ribbon, but {\em every} knot is ribbon. 

The point of confusion lies in the precise nature of the Douglas-Rad{\'o} solution.  They showed that any rectifiable curve in $\rr^n$ (in particular, every knot in $S^3\subset \rr^4$)  arises as the boundary of a map $u \co D^2\rightarrow \rr^n$ which is absolutely area-minimizing, and hence geometrically minimal, but may have two types of singularities.  The first type are self-intersections; that is, the map is at best an immersion. The second type of singularity, called {\em branch points}, are  the image of points  $p\in D^2$ at which the derivative  vanishes: $du_p=0$.   It is this latter type of singularity which ties the ribbon-slice problem to the concordance problem of algebraic knots.

The key point is that the structure  of branch points of minimal surfaces in $\rr^4$ can be understood in terms very similar to those in the classical theory of plane curve singularities.  Indeed, Micallef and White~\cite{MW95} showed that the link of a branch point of an area-minimizing map $u \co D^2\rightarrow \rr^4$ is equivalent, up to a possibly  orientation-reversing diffeomorphism of $S^3=\partial D^4$, to an algebraic link.  Thus the boundary of a small neighborhood of the pre-image of a branch point, $N_\epsilon(p)\subset D^2$,  is mapped into the three-sphere $\partial N_\delta(u(p))$ as the link type of an algebraic link or its mirror image. (Here $N_\delta(u(p))$ is a neighborhood in $\rr^4$ of the branch point.)

In the case of a topologically embedded minimal disk we can group any branch points  together to obtain a  ribbon concordance (that is, a ribbon cylinder) between the boundary of the disk and a connected sum of knots, $K_1\#...\#K_n$,  where each $K_i$ is either an algebraic knot or the mirror image of an algebraic knot.  If we could show that the only such sums which are slice are of the form $K\#-K$, then we could replace the neighborhood of the collection of branch points in the minimal disk with the obvious ribbon disk bounded by such connected sums.   This replacement results in a ribbon disk for the original knot.   

Thus the ribbon slice problem could be proved in the affirmative with two (difficult) steps.  
\begin{enumerate} 
\item Show that any slice knot type can be realized as the boundary of a topologically embedded minimal disk, but possibly with branch points. 
\item Classify concordance relations between algebraic links and show that the only such relations are of the form $K\#-K$. Note that if any other relations exist, they cannot be realized by ribbon disks, by Miyazaki~\cite{Miyazaki94}.  In this case the answer to the ribbon-slice question would be in the negative.
\end{enumerate}

 The first step seems more difficult  due to the nature of known results for Plateau-type problems, which  indicate that self-intersections of minimal maps are more likely to occur generically than branch points. Hence it may be too much to hope that one can find the required minimal embeddings through the current techniques.  It should be pointed out, however, that Hass~\cite{Hass} showed that any ribbon knot has an isotopy representative bounding an embedded minimal disk without branch points.  Thus if the answer to the ribbon-slice question is affirmative, then the strategy works without the second step.  It is our hope, however, that importing the concordance result could provide useful flexibility in the geometric approach to the problem.

\end{document}